%meta.tex Plain TeX file by Eric Rowland and Doron Zeilberger (xx pages)

\input amssym.tex

%begin macros
\baselineskip=14pt
\parskip=10pt

\font\eighttt=cmtt8
\magnification=\magstephalf
\def\N{{\Bbb N}}
\def\Z{{\Bbb Z}}
\def\1{{\overline{1}}}
\def\2{{\overline{2}}}
\parindent=0pt
\overfullrule=0in

\def\frac#1#2{{#1 \over #2}}

%\def\RightHead{\centerline{
%Title
%}}
%\def\LeftHead{ \centerline{Doron Zeilberger}}
%end macros

\centerline
{
\bf A Case Study in Meta-AUTOMATION:  
}
\centerline
{\bf
AUTOMATIC Generation of  Congruence AUTOMATA For Combinatorial Sequences
}
\bigskip
\centerline
{
\it Eric ROWLAND and Doron ZEILBERGER
}

\quad \quad \quad \quad \quad

{\bf Abstract}: In this paper, that may be considered a sequel to a recent article
by Eric Rowland and Reem Yassawi,
we present yet another approach for the automatic generation of automata
(and an extension that we call Congruence Linear Schemes) for the fast (log-time) determination
of congruence properties, modulo small (and not so small!)\ prime powers,  
for a wide class of combinatorial sequences.
Even more interesting than the new results that could be obtained, is the illustrated
{\it methodology}, that of designing `meta-algorithms' that enable the computer to
develop algorithms, that it (or another computer) can then proceed to use to actually
prove (potentially!)\ infinitely many new results.
This paper is accompanied by a Maple package, {\tt AutoSquared}, and numerous sample
input and output files, that readers can use as templates for generating
their own, thereby proving many new `theorems' about congruence properties
of many famous (and, of course, obscure) combinatorial sequences.

{\bf Very Important}: This article is accompanied by the general Maple package

{\eighttt http://www.math.rutgers.edu/\~{}zeilberg/tokhniot/AutoSquared} \quad ,

and several other specific ones, and numerous input and output 
files that are obtainable, by {\it one} click, from the webpage (``front'') of this article

{\tt http://www.math.rutgers.edu/\~{}zeilberg/mamarim/mamarimhtml/meta.html} \quad \quad .

They could (and should!)\ be used as templates for generating as many input files
that the human would care to type, and that the computer would agree to run.

{\bf Prologue: What are the Last Three (Decimal) Digits of the Googol-th Catalan, Motzkin, and Delannoy Numbers?}

We will {\it never} know {\it all} the (decimal) digits of the Googol-th terms of the famous 
{\it Catalan}, {\it Motzkin}, and (Central) {\it Delannoy} sequences
[{\tt http://oeis.org/A000108},  {\tt http://oeis.org/A001006}, and {\tt http://oeis.org/A001850}, respectively],
if nothing else because our computers are not big enough to store them!

But thanks to the Maple packages accompanying this article, we know {\it for sure} that the {\bf last three digits} are
$000$, $187$, and $281$ respectively. These packages can compute in {\it logarithmic time}
(i.e.\ linear in the number of digits of the input) the values of the $n$-th term modulo many different $m$
(but alas, not too big!). These fast algorithms were generated by a {\it meta-algorithm} implemented
in the main Maple package, {\tt AutoSquared}.

{\bf Fast Exponentiation}

E-commerce is possible (via RSA) thanks to the fact that it is very easy (for computers!)\ to compute
$$
a^n \bmod{m}  \quad  \quad ,
$$
for $a$ and $m$ several-hundred-digits long, and large $n$. Reminding you that $a^n$ is
shorthand for the {\it sequence}, let's call it $x_n$ defined by the {\bf linear recurrence equation}
with {\bf constant} coefficients, of {\bf order} one:
$$
x_n -a x_{n-1}=0 \quad , \quad x_0=1 \quad .
$$
In order to compute $x_{10^{100}} \bmod{m}$, you don't compute all the $10^{100}$ previous terms, but
use the {\bf implied} recurrences
$$
x_{2n} = x_n^2 \bmod{m} \quad , \quad
x_{2n+1} = a x_n^2 \bmod{m}  \quad \quad .
$$
This takes only $\log_2 10^{100}$ operations!

What about sequences defined by higher-order recurrences, but still with {\it constant} coefficients?
For example, what are the last three decimal digits of the googol-th Fibonacci number, $F_{10^{100}}$?
You would get the answer, $875$, in {\tt 0.008} seconds!

All you need is type

{\tt Fnm(10**100, 1000); }  \qquad ,

once you typed (or copied-and-pasted) the following short code into a Maple session:

{\tt
Fnm:=proc(n, m) option remember; \hfill\break
    if n = 1 or n = 2 then 1 \hfill\break
    elif n mod 2 = 0 then Fnm(1/2*n, m)*(Fnm(1/2*n + 1, m) + Fnm(1/2*n - 1, m)) mod m \hfill\break
    else Fnm(1/2*n - 1/2, m)**2 + Fnm(1/2*n + 1/2, m)**2 mod m \hfill\break
    fi: \hfill\break
end: \hfill\break
}

It implements the (nonlinear) recurrence scheme
$$
F_{2n}=F_n (F_{n-1}+F_{n+1})  \quad , \quad  F_{2n+1}=F_{n}^2 + F_{n+1}^2  \quad , \quad  F_1=1 \quad , \quad F_2=1 \quad ,
$$
and of course takes it modulo $m$ at every step.

Another way is to take the $(1,2)$ entry of the matrix
$$
\pmatrix{ 1& 1 \cr 1 & 0 }^{10^{100}} \bmod{1000} \quad,
$$
and use the `iterated-squaring' trick applied to  matrix (rather than scalar) exponentiation.

Both these simple methods are applicable for the fast (linear-in-bit-size) computation of
the terms, modulo {\it any} $m$,  of {\it any}  integer sequence defined in terms
of a {\it linear recurrence equation} with {\bf constant} coefficients (aka $C$-finite integer sequences).

But what about sequences that are defined via {\it linear recurrence equations with} {\bf polynomial} coefficients,
aka {\it P-recursive} sequences, aka {\it holonomic} sequences?

In a beautiful and deep paper ([KKM]), dedicated to one of us (DZ) on the occasion of his $60^{th}$ birthday,
Manuel Kauers, Christian Krattenthaler, and Thomas M\"uller developed
a deep and ingenious group-theoretical method for the determination of holonomic sequences modulo powers of $2$.
This has been extended to powers of $3$ in [KM1], and further developed in [KM2].

An important subclass of the class of holonomic integer sequences is the class of integer sequences whose (ordinary) generating
function, let's call it $f(x)$, satisfies an algebraic equation of the form $P(f(x),x)=0$ where
$P$ is a polynomial of two variables with integer coefficients. For this class, and an even
wider class, the sequences arising from the {\it diagonals} of rational functions of several variables,
Rowland and Yassawi ([RY]) developed a very {\it general} method for computing {\bf finite automata} 
for the fast computation (once the automaton is found, of course) of the congruence behavior
modulo prime powers. Of course, as the primes and/or their powers get larger, the automata get
larger too, but if the automaton  is precomputed {\it once and for all} (and saved!), it is logarithmic time
(i.e.\ linear in the bit-size). Of course, the {\it implied constant} in the $O(\log n)$ computation
times gets larger with the moduli.

{\bf History}

Many papers, in the past, proved isolated results about congruence  properties for {\it specific} sequences
and for {\it specific} moduli. We refer the reader to [RY] for many  references, that we will not 
repeat here.

{\bf The Present Method: Using Constant Terms}

Most (perhaps all) of the combinatorial sequences treated in [RY] can be written in the form
$$
a_n := ConstantTermOf \, [ \, P(x)^n Q(x) \,] \quad ,
$$
where both $P(x)$ and $Q(x)$ are {\it Laurent polynomials} with integer coefficients, where
$x$ is either a single variable or a multi-variable $x=(x_1, \dots, x_m)$, and ConstantTermOf means
``coefficient of $x^0 \,\,$'', or ``the coefficient of $x_1^0 \cdots x_m^0 \,\,$''.

For example, the arguably second-most famous combinatorial sequence (after the Fibonacci sequence), 
is the sequence of the Catalan Numbers ({\tt http://oeis.org/A000108}), that may be defined by
$$
C_n :=ConstantTermOf \,\,[\,\, (\frac{1}{x}+2+x)^n (1-x)\,\,] \quad .
$$
Not as famous, but also popular, are the {\bf Motzkin} numbers ({\tt http://oeis.org/A001006}), that may be defined by
$$
M_n :=ConstantTermOf \,\,[ \,\, (\frac{1}{x}+1+x)^n (1-x^2) \,\,] \quad ,
$$
and also fairly famous are the {\bf Central Delannoy Numbers} ({\tt http://oeis.org/A001850}), that may be defined
by
$$
D_n :=ConstantTermOf \,\, [ \,\, (\frac{1}{x}+3+2x)^n \,\, ] \quad .
$$
So far, we got away with a single variable. 

Another celebrated sequence is the sequence of {\it Ap\'ery Numbers}, that
were famously used by 64-year-old Roger Ap\'ery (in 1978) to prove the
irrationality of $\zeta(3)$. These are defined in terms of
a {\it binomial coefficient sum}
$$
A(n):= \sum_{k=0}^{n} {{n} \choose {k}}^2 {{n+k} \choose {k}}^2 \quad .
$$
These may be equivalently defined (see below) as
$$
A(n):=
ConstantTermOf \left [
\left (
{\frac { \left( 1+x_{{1}} \right)  \left( 1+x_{{2}} \right)  \left( 1+x_{{3}} \right)  
\left(1+x_{{2}}+x_{{3}}+x_{{2}}x_{{3}}+  x_{{1}}x_{{2}}x_{{3}}
 \right) }{x_{{1}}x_{{2}}x_{{3}}}} \right )^n \right ] \quad .
$$

{\bf How to convert ANY Binomial Coefficient Sum into a Constant Term Expression?}

Before describing our new method, let us indicate how any binomial coefficient sum
of the form
$$
A(n)= \sum_{k=0}^{n}  {{n} \choose {k}}  g^k \prod_{i=1}^{m}  {{a_i n +b_i k + c_i} \choose {d_i n +e_i k + f_i}} \quad,
$$
where all the $a_i,b_i,c_i,d_i,e_i,f_i$ and $g$ are {\bf integers}, can be made into a constant term expression.
(This is essentially Georgy Petrovich EGORYCHEV's celebrated {\it method of coefficients}).
We introduce $m$ variables $x_1, \dots, x_m$ and use the fact, that, {\it by definition}
$$
 {{a_i n +b_i k + c_i} \choose {d_i n +e_i k + f_i}}=
ConstantTermOf_{x_i} \, \left [ \, \frac{(1+x_i)^ {a_i n +b_i k + c_i}}{x_i^{d_i n +e_i k + f_i}} \, \right ]
$$
Hence
$$
A(n)= \sum_{k=0}^{n}  {{n} \choose {k}} g^k \prod_{i=1}^{m}  {{a_i n +b_i k + c_i} \choose {d_i n +e_i k + f_i}} 
$$
$$
=
\sum_{k=0}^{n}  {{n} \choose {k}} g^k
\prod_{i=1}^{m}  ConstantTermOf_{x_i} \,  [ \, \frac{(1+x_i)^ {a_i n +b_i k + c_i}}{x_i^{d_i n +e_i k + f_i}} \, ] 
$$
$$
=
ConstantTermOf_{x_1, \dots , x_m} \, \left [ \,
\sum_{k=0}^{n}  {{n} \choose {k}} g^k
\prod_{i=1}^{m}  \frac{(1+x_i)^ {a_i n +b_i k + c_i}}{x_i^{d_i n +e_i k + f_i}} 
\, \right ]
$$
$$
=
ConstantTermOf_{x_1, \dots , x_m} \, \left [ \,
\left ( \prod_{i=1}^{m}  \frac{(1+x_i)^ {a_i n + c_i}}{x_i^{d_i n + f_i}} \right ) 
\sum_{k=0}^{n}  {{n} \choose {k}} g^k
\prod_{i=1}^{m} \left ( \frac{(1+x_i)^ {b_i k}}{x_i^{e_i k}} \right ) 
\, \right ]
$$
$$
=
ConstantTermOf_{x_1, \dots , x_m} \, \left [ \,
 \left ( \, \prod_{i=1}^{m} \frac{(1+x_i)^ {a_i n + c_i}}{x_i^{d_i n + f_i}} \, \right ) 
\sum_{k=0}^{n}  {{n} \choose {k}} g^k
\prod_{i=1}^{m} \left ( \frac{(1+x_i)^ {b_i }}{x_i^{e_i }} \right )^k 
\, \right ]
$$
$$
=
ConstantTermOf_{x_1, \dots , x_m} \, \left [ \,
\prod_{i=1}^{m} \left ( \frac{(1+x_i)^ {a_i n + c_i}}{x_i^{d_i n + f_i}} \right ) 
 \left ( 1+g \prod_{i=1}^{m}  \frac{(1+x_i)^ {b_i }}{x_i^{e_i }} \right )^n
\, \right ] \quad .
$$
$$
=
ConstantTermOf_{x_1, \dots , x_m} \, \left [ \,
\prod_{i=1}^{m} \left ( \frac{(1+x_i)^ { c_i}}{x_i^{f_i}} \right ) 
\left ( \frac{(1+x_i)^ {a_i }}{x_i^{d_i}} \right )^n
 \left ( 1+g \prod_{i=1}^{m}  \frac{(1+x_i)^ {b_i }}{x_i^{e_i }} \right )^n
\, \right ] \quad .
$$
$$
=
ConstantTermOf_{x_1, \dots , x_m} \, \left [ \,
\prod_{i=1}^{m}  \frac{(1+x_i)^{ c_i}}{x_i^{f_i}} 
\left ( \left ( \prod_{i=1}^{m} \frac{(1+x_i)^{a_i}}{x_i^{d_i}}  \right )
\left ( 1+g \prod_{i=1}^{m}  \frac{(1+x_i)^{b_i }}{x_i^{e_i }} \right ) \right )^n
\, \right ] \quad .
$$
This is implemented in procedure {\tt BinToCT(L,x,a)} in our Maple package {\tt AutoSquared}. For example,
we got the above constant-term rendition of the Ap\'ery numbers by typing:

{\tt  BinToCT([ [[1,0,0],[0,1,0]], [[1,1,0],[0,1,0]]\$2 ],x,1); } .

{\bf Illustrating the Constant Term Approach In Terms of the Simplest-Not-Entirely-Trivial Example}

Recall from above that the Catalan numbers may be defined by the constant-term formula
$$
C_n :=ConstantTermOf \,\,[\,\, (\frac{1}{x}+2+x)^n (1-x)\,\,] \quad .
$$
We are interested in the mod $2$ behavior of $C_n$, in other words we want to have
a quick way of computing $C_n$ modulo $2$. So let's define
$$
A_1(n) := C_n \bmod{2} \quad .
$$
Using the above formula for $C_n$, and taking it modulo $2$, we have:
$$
A_1(n) :=ConstantTermOf \,\,[\,\,  (1+x) (\frac{1}{x}+x)^n\,\,] \quad .
$$
We will try to find a constant-term expression for $A_1(2n)$ .
$$
A_1(2n) =ConstantTermOf \,\,[\,\,  (1+x) (\frac{1}{x}+x)^{2n}\,\,]  \bmod 2
=ConstantTermOf \,\,[\,\,  (1+x) ((\frac{1}{x}+x)^{2})^{n}\,\,] \bmod 2
$$
$$
=ConstantTermOf \,\,[\,\,  (1+x) (\frac{1}{x^2}+2+x^{2})^{n}\,\,] \bmod 2
=ConstantTermOf \,\,[\,\,  (1+x) (\frac{1}{x^2}+x^{2})^{n}\,\,] \bmod 2
$$
But
$$
ConstantTermOf \,\,[\,\,  (1+x) (\frac{1}{x^2}+x^{2})^{n}\,\,] 
\, = \, ConstantTermOf \,\,[\,\,  1 \cdot (\frac{1}{x^2}+x^{2})^{n}\,\,] \quad ,
$$
since, obviously,
$$
ConstantTermOf \,\,[\,\,  x \cdot (\frac{1}{x^2}+x^{2})^{n}\,\,] =0 \quad .
$$
Since  the constant-termand of 
$$
ConstantTermOf \,\,[\,\,  1 \cdot (\frac{1}{x^2}+x^{2})^{n}\,\,] \quad ,
$$

{\bf only} depends on $x^2$, we can replace $x^2$ by $x$, implying that
$$
A_1(2n) =ConstantTermOf \,\,[\,\,  (\frac{1}{x}+x)^{n}\,\,]  \bmod 2 \quad .
$$
This {\bf forces} us to {\bf put up} with {\bf a new kid on the block}, let's call it
$A_2(n)$:
$$
A_2(n) :=ConstantTermOf \,\,[\,\,   (\frac{1}{x}+x)^n\,\,] \bmod  2\quad ,
$$
and we got the {\it recurrence}
$$
A_1(2n)=A_2(n) \quad .
$$
We will handle $A_2(n)$ in due course, but first let's
consider $A_1(2n+1)$.

We have
$$
A_1(2n+1) =ConstantTermOf \,\,[\,\,  (1+x) (\frac{1}{x}+x)^{2n+1}\,\,]  \bmod 2
$$
$$
=ConstantTermOf \,\,[\,\,  (1+x)(\frac{1}{x}+x) ((\frac{1}{x}+x)^{2})^{n}\,\,] \bmod 2
$$
$$
=ConstantTermOf \,\,[\,\,  (\frac{1}{x}+x+1+x^2) (\frac{1}{x^2}+2+x^{2})^{n}\,\,] \bmod 2
$$
$$
=ConstantTermOf \,\,[\,\,    (\frac{1}{x}+x+1+x^2) (\frac{1}{x^2}+x^{2})^{n}\,\,] \bmod 2 \quad .
$$
But
$$
ConstantTermOf \,\,[\,\,    (\frac{1}{x}+x+1+x^2) (\frac{1}{x^2}+x^{2})^{n}\,\,] 
\, = \, ConstantTermOf \,\,[\,\,  (1 +x^2) \cdot (\frac{1}{x^2}+x^{2})^{n}\,\,] \quad ,
$$
since, obviously
$$
ConstantTermOf \,\,[\,\,  (\frac{1}{x} + x) \cdot (\frac{1}{x^2}+x^{2})^{n}\,\,] =0 \quad .
$$
Since  the constant-termand of 
$$
ConstantTermOf \,\,[\,\,  (1+x^2) \cdot (\frac{1}{x^2}+x^{2})^{n}\,\,] \quad ,
$$
{\bf only} depends on $x^2$, we can replace $x^2$ by $x$, implying that
$$
A_1(2n+1) =ConstantTermOf \,\,[\,\,  (1+x) (\frac{1}{x}+x)^{n}\,\,]  \bmod 2 \quad .
$$
But this looks familiar! It is good-old $A_1(n)$, so we have established, so far, that
$$
A_1(2n)=A_2(n) \quad , \quad A_1(2n+1)=A_1(n) \quad .
$$
But in order to establish a {\bf recurrence scheme}, we need to handle $A_2(n)$. A priori, this
may force us to introduce yet more discrete functions, and that would be OK, as long as we would
{\it finally} stop, after {\it finitely} many steps, getting a scheme with {\bf finitely} many
discrete functions, that would enable the fast  (logarithmic time) computation of 
our initial function $A_1(n)$. We will see that this would {\it always} be the case, 
no matter how complicated $P(x)$ and $Q(x)$ are (and even with many variables). Alas,
as $P(x)$ gets more complicated, the `finite' gets bigger and bigger, so eventually the
`logarithmic time' in $n$ would be impractical, since the {\it implied constant} would be eeeeeeeeeeeeenormous.

But in this {\bf toy example}, don't worry! The `finitely many  discrete functions', is only two!
As we will shortly see, all we need is $A_2(n)$, in addition to $A_1(n)$.

Recall that
$$
A_2(n) :=ConstantTermOf \,\,[\,\,   (\frac{1}{x}+x)^n\,\,] \bmod 2 \quad .
$$
Let's try to find a constant-term expression for $A_2(2n)$ .
$$
A_2(2n) =ConstantTermOf \,\,[\,\,  (\frac{1}{x}+x)^{2n}\,\,]  \bmod 2
=ConstantTermOf \,\,[\,\,  ( (\frac{1}{x}+x)^{2})^{n}\,\,] \bmod 2
$$
$$
=ConstantTermOf \,\,[\,\,  (\frac{1}{x^2}+2+x^{2})^{n}\,\,] \bmod 2
=ConstantTermOf \,\,[\,\,   (\frac{1}{x^2}+x^{2})^{n}\,\,] \bmod 2
$$

Since the constant-termand {\bf only} depends on $x^2$, we can replace $x^2$ by $x$,
implying that
$$
A_2(2n) =ConstantTermOf \,\,[\,\,  (\frac{1}{x}+x)^{n}\,\,]  \bmod 2 \quad .
$$
But that's exactly $A_2(n)$, so we have found out that
$$
A_2(2n)=A_2(n) \quad .
$$

What about $A_2(2n+1)$? Here goes:
$$
A_2(2n+1) =ConstantTermOf \,\,[\,\,  (\frac{1}{x}+x)^{2n+1}\,\,]  \bmod 2
=ConstantTermOf \,\,[\,\, (\frac{1}{x}+x) ((\frac{1}{x}+x)^{2})^{n}\,\,] \bmod 2
$$
$$
=ConstantTermOf \,\,[\,\,   (\frac{1}{x}+x) (\frac{1}{x^2}+2+x^{2})^{n}\,\,] \bmod 2
=ConstantTermOf \,\,[\,\,   (\frac{1}{x}+x) (\frac{1}{x^2}+x^{2})^{n}\,\,] \bmod 2 \quad .
$$
But the constant-termand now only has {\bf odd} powers, so the coefficient of $x^0$, alias the
{\it constant term}, is $0$. We have just established, the {\bf fast recurrence scheme}:
$$
A_1(2n)=A_2(n) \quad , \quad A_1(2n+1)=A_1(n) \quad ,
$$
$$
A_2(2n)=A_2(n) \quad , \quad A_2(2n+1)=0 \quad ,
$$
subject to the {\bf initial conditions}
$$
A_1(0)=1 \quad, \quad A_2(0)=1 \quad.
$$

[The above human-generated scheme can be also done ({\bf much faster}) by the Maple package \hfill\break
{\tt AutoSquared}.
Having downloaded  {\tt http://www.math.rutgers.edu/\~{}zeilberg/tokhniot/AutoSquared} into
your computer, that has Maple installed, you stay in the same directory, and you type:

{\tt  read AutoSquared: CA([1/x+2+x,1-x],x,2,1,2)[1];} \qquad ,

and you would get (in 0 seconds!), the output 

{\tt [[[2, 1], [2, 0]], [1, 1]]} \qquad ,

which is our package's way of encoding the above `scheme'.

Another way of describing the scheme is via the {\it binary} representation of $n$ (for some $k \geq 1$)
$$
n=\sum_{i=1}^{k} \alpha_i 2^{k-i} \quad ,
$$
where $\alpha_i \in \{0,1\}$, $\alpha_1=1$, and it is abbreviated, in the {\bf positional} notation, as a {\bf word},
of length $k$, in the {\bf alphabet} $\{0,1\}$
$$
\alpha_1 \cdots \alpha_k \quad .
$$
Phrased in terms of such `words', the above scheme can be written, 
(where $w$  is {\it any} word in the alphabet $\{0,1\}$)
$$
A_1(w0)=A_2(w) \quad , \quad A_1(w1)=A_1(w) \quad .
$$
$$
A_2(w0)=A_2(w) \quad , \quad A_2(w1)=0 \quad ,
$$
subject to the {\it initial conditions} (here $\phi$ is the empty word):
$$
A_1(\phi)=1 \quad, \quad A_2(\phi)=1 \quad .
$$

Let's revert to {\it post-fix} notation for representing functions, and omit the parentheses, i.e.\
write $wA_1$ instead of $A_1(w)$ and $wA_2$ instead of $A_2(w)$. This will not cause any ambiguity,
since the {\bf alphabet of function names}, $\{ A_1, A_2 \}$ is disjoint from the {\bf alphabet of letters},
$\{0,1\}$.
The above scheme becomes
$$
w0A_1=wA_2 \quad , \quad w1A_1=wA_1 \quad
$$
$$
w0A_2=wA_2 \quad , \quad w1A_2=0 \quad ,
$$
subject to the {\it initial conditions}
$$
\phi A_1=1 \quad, \quad \phi A_2=1 \quad.
$$

Let's try 
to find $A_1(30)$, alias,  $A_1(11110_2)$, alias, with our new convention,
$11110A_1$. We get in {\bf two} steps 
$$
11110A_1=1111A_2=0 \quad .
$$
This only took two steps due to a premature exit to an output gate. The default
number of steps is the length of the word, that keeps traveling until
it becomes the empty word, and then it is {\bf forced} to move to an output gate.

It is readily seen that if the input word has a zero in it, the output would be $0$. Hence the
only words that output $1$ are those given by the {\bf regular expression} 
$$
1^{*} \quad .
$$
Equivalently, the only integers $n$ for which the Catalan number $C_n$ is odd are those of the form
$n=2^k-1$ for $k=0,1,2, \dots$.

The words in the alphabet $\{0,1\}$ that output $0$ (i.e.\ those words that have at least one $0$ in their
binary representation) are the complement `language', whose
regular expression rendition is
$$
1\, \{0,1\}^{*} \, 0 \, \{0,1\}^{*} \quad .
$$

What we have here is a {\bf finite automaton with output}. The set of {\bf states} is $\{A_1,A_2\}$ while the
{\bf alphabet} is the set $\{0,1\}$. There are $2$ 
directed edges coming out of each state, one for each letter of the alphabet,
leading to another (possibly the same) state, or possibly to an output gate (in our case always $0$,
via `exit edges' that prematurely end the journey.
You have a {\it starting state} (in this example, $A_1$) and an input word, and you travel along the
automaton, according to the current state and the current rightmost letter, until you run out of letters,
i.e.\ have the empty word, or wind-up in the output $0$ prematurely, since some states have edges that
lead directly to $0$. (In our example when you are at state $A_2$ and the rightmost letter is $1$
you immediately output $0$.)

Yet another way of describing it is via a {\bf type-three grammar} (aka {\bf regular grammar}) in the famous {\bf Chomsky hierarchy}
(see e.g.\ [R]).
For each possible output (in this example, $0$ and $1$, {\bf NOT TO BE CONFUSED WITH THE LETTERS OF THE ALPHABET}), there is a regular grammar
describing the language (set of words) that yield that output. 

In this example, the set of {\bf non-terminal symbols} is $\{A_1,A_2\}$ and the set of  {\bf terminal symbols} is
$\{0,1\}$. For a grammar for the language yielding $1$ (i.e.\ the binary representations of the integers $n$
for which $C_n$ is odd) the non-terminal symbol $A_2$ is not needed (is superfluous), and the grammar is extremely simple
$$
A_1 \rightarrow \phi \quad , \quad A_1 \rightarrow 1 A_1 \quad .
$$
We leave it to the interested reader to write down the only slightly more complicated grammar for the
language of binary representations of integers $n$ for which $C_n$ is even.

It is well-known that the notions of {\it finite automata}, {\it regular expressions}, and {\it regular grammars}
are equivalent (as far as the generated languages), and there are easy algorithms for going between them.

These are all very nice, but for the present formulation, it is more convenient {\bf not} to write the input integers
$n$ in base $2$ (or more generally, base $p$, if the desired modulus is a power of a prime $p$), but
stick to integers (as inputs). Let's make the following formal definition.

{\bf Definition}: Let $\N$ be the set of non-negative integers, let $p$ be a positive integer,
and let $E$ be any set.
An  {\it automatic $p$-scheme} for a function $f:\N \rightarrow E$ is
a set of {\it finitely} many (say $r$) auxiliary functions $A_1(n), \dots, A_r(n)$, where
$f(n)=A_1(n)$ and there is a function
$$
\sigma: \{0, \dots, p-1\} \times \{1, \dots, r\} \rightarrow  \{1, \dots, r\}  \quad ,
$$
such that, for each $1 \leq i \leq r$ and $0 \leq \alpha  \leq p-1$, we have the recurrence
$$
A_i(pn+\alpha)= A_{\sigma(\alpha,i)}(n)  \quad .
$$
We also have {\bf initial conditions} 
$$
A_i(0)=a_i \quad ,
$$
for some $a_i \in E \quad, 1 \leq i \leq r$.

{\bf Note}: In the application to schemes for congruence properties of combinatorial sequences 
modulo prime powers $p^a$, treated in the present article, $p$ will always be a prime, and the output set, $A$, would be
$$
\{ 0,1, \dots, p^a-1 \} \quad .
$$

{\bf Teaching the Computer How to Create Automatic $p$-schemes}

All the tricks described above, in excruciating detail, for finding the scheme for determining the
mod $2$ behavior of the Catalan numbers 
$$
C_n :=ConstantTermOf \,\,[\,\, (\frac{1}{x}+2+x)^n (1-x)\,\,] \quad ,
$$
can be taught to the computer (in our case using the symbolic programming language Maple), to find
{\it without human touch}, an automatic $p$-scheme for determining the
mod $p^a$ behavior, for {\bf any} prime $p$, and {\bf any} power $a$, for {\bf any}
combinatorial sequence defined by
$$
A(n) :=ConstantTermOf \, [ \, P(x_1, \dots, x_m)^n Q(x_1, \dots, x_n) \,] \bmod p^a \quad ,
$$
for {\bf any} polynomials with integer coefficients, $P(x_1, \dots, x_m)$ and $Q(x_1, \dots, x_m)$,
for {\bf any} number of variables.

We will associate $A(n)$ with the pair $[P,Q]$.

We first rename $A(n)$, $A_1(n)$, and $[P,Q]$, $[P_1,Q_1]$.
We then try to find constant-term expressions for $A_1(np), A_1(np+1), \dots, A_1(np+p-1)$.
After using the multinomial theorem and doing it mod $p^a$, we would get, e.g.,
$$
A_1(pn) =ConstantTermOf \, [ \, P_1(x_1, \dots, x_m)^{np} Q_1(x_1, \dots, x_{m}) \,] \bmod p^a 
$$
$$
=ConstantTermOf \, [ \, (P_1(x_1, \dots, x_m)^p)^n Q_1(x_1, \dots, x_{m}) \,] \bmod p^a \quad ,
$$
that after simplification (expanding, taking modulo $p^a$, and, if applicable, replacing $x^p$ by $x$)
will force us to put up with a brand-new discrete function, let's call it $A_2(n)$, given by
$$
A_2(n)=ConstantTermOf \, [ \, P_2(x_1, \dots, x_m)^n Q_2(x_1, \dots, x_{m}) \,] \bmod p^a \quad ,
$$
So $A_2$ corresponds to a brand-new pair $[P_2,Q_2]$. We do likewise for $A_1(pn+1)$, all the way to $A_1(pn+p-1)$,
getting (at the beginning) new pairs.
Then we do the same for $A_2(pn)$ through $A_2(pn+p-1)$.
After awhile, by the {\bf pigeonhole principle}, we will get old friends, and eventually there won't be any `new guys', and we get a {\bf finite} (alas, often very large!)\
automatic $p$-scheme.
The proof is as follows.
If $P(x)$ is a Laurent polynomial in $x_1^p, \dots, x_m^p$, let $\Lambda(P(x))$ denote the Laurent polynomial obtained by replacing each $x_j^p$ by $x_j$.
Since $\Lambda$ commutes with raising to the $p$th power, the first component of each pair $[P_i, Q_i]$ after $a$ iterations is $\Lambda^k(P(x)^{p^a})$ for some $k \geq 0$.
The only terms of $P(x)^{p^a}$ whose coefficients are non-zero modulo $p^a$ are those in which the exponent of each $x_j$ is a multiple of $p$; therefore $k \geq 1$.
It is not too difficult to see (for example, using Proposition~1.9 in [RY]) that
$$
	\Lambda\left(P(x)^{p^a}\right)
	\equiv P(x)^{p^{a-1}} \pmod {p^a}.
$$
From this it follows that
$$
	\Lambda^k\left(P(x)^{p^a}\right)
	\equiv \Lambda^{k-1}\left(P(x)^{p^{a-1}}\right) \pmod {p^a}.
$$
On the next iteration, we raise this polynomial to the $p$th power and apply $\Lambda$; this gives
$$
	\Lambda\left(\left(\Lambda^{k-1}(P(x)^{p^{a-1}})\right)^p\right) \bmod p^a
	= \Lambda^k\left(P(x)^{p^a}\right) \bmod p^a
	= \Lambda^{k-1}\left(P(x)^{p^{a-1}}\right) \bmod p^a,
$$
so the first component of $[P_i, Q_i]$ stays the same after $a$ iterations.
There are only finitely many possibilities for the second component as well, since after the first component stabilizes then we can apply $\Lambda$ to both $P$ and (after deleting some terms) $Q$ at each iteration, and this puts bounds on the degree and low-degree of $Q$.

All of this is implemented in {\tt AutoSquared} by procedure {\tt CA} for single-variable 
polynomials $P$ and $Q$ and
by procedure {\tt CAmul} for multivariate $P$ and $Q$ (of course, {\tt CAmul} can handle also
a single variable, but we kept {\tt CA} both for old-time-sake and because it may be a bit faster for this special case).

The syntax is

{\tt CA(Z,x,p,a,K): } \qquad ,

where {\tt Z} is a pair of single-variable functions $[P,Q]$, 
in the variable {\tt x}, {\tt p} is a prime, {\tt a} is a positive integer, and {\tt K} is 
a (usually large) positive integer,
stating the maximum number of `states' (auxiliary functions) that you are willing to put up with.
(It returns FAIL if the number of states exceeds {\tt K}.)

For example, to get an automatic $2$-scheme for the Motzkin numbers, modulo $2$, 
(if you are willing to tolerate a scheme with up to $30$ members),

you type: \quad {\tt gu:=CA([1/x+1+x,1-x**2],x,2,1,30): } \quad .

The output (that we call {\tt gu}) has two parts. The second part, {\tt gu[2]}, that is {\bf not} needed
for the application for the fast determination of the sequence modulo $2$ (and in general modulo $p^a$)
consists in the `definition', in terms of constant term expressions $A_i(n):=ConstantTerm[P_i(x)^n Q_i(x)]$,
of the various auxiliary functions. So, in this example, {\tt gu[2]} is

{\tt [[1/x+1+x, 1+x**2], [1/x+1+x, 1+x], [1/x+1+x, 1], [1/x+1+x, x]]} \quad,

meaning that 
$$
A_1(n)=ConstantTerm[(1/x+1+x)^n(1+x^2)] \quad , \quad A_2(n)=ConstantTerm[(1/x+1+x)^n(1+x)]
$$
$$
A_3(n)=ConstantTerm[ (1/x+1+x)^n ] \quad , \quad A_4(n)=ConstantTerm[(1/x+1+x)^n \cdot x)] \quad .
$$
The more interesting part, the one needed for the actual fast computation, is {\tt gu[1]}.

Typing : \quad {\tt lprint(gu[1])} \quad in the same Maple session, gives

{\tt [[[2, 2], [3, 4], [3, 3], [0, 2]], [1, 1, 1, 0]]} \quad ,

that in {\it humanese} means the $2$-scheme
$$
A_1(2n)=A_2(n) \quad , \quad  A_1(2n+1)=A_2(n) \quad , 
$$
$$
A_2(2n)=A_3(n) \quad , \quad  A_2(2n+1)=A_4(n) \quad , 
$$
$$
A_3(2n)=A_3(n) \quad , \quad  A_3(2n+1)=A_3(n) \quad , 
$$
$$
A_4(2n)=0 \quad , \quad  A_4(2n+1)=A_2(n) \quad .
$$
The initial conditions are
$$
A_1(0)=1 \quad , \quad A_2(0)=1 \quad , \quad A_3(0)=1 \quad , \quad A_4(0)=0 \quad .
$$

Moving right along, to get an automatic $2$-scheme for the Motzkin numbers mod $4$
(let's tolerate from now on systems up to $10000$ states):

{\tt gu:=CA([1/x+1+x,1-x**2],x,2,2,10000): } \quad ,

getting (by typing {\tt nops(gu[2])} (or {\tt nops(gu[1][1])})) a scheme with $24$ states.

To get an automatic $2$-scheme for the Motzkin numbers mod $8$ (still with $\leq 10000$ states, if possible), you type

{\tt gu:=CA([1/x+1+x,1-x**2],x,2,3,10000): } \quad ,

getting a certain scheme with $128$ states.

For mod $16$, we type

{\tt gu:=CA([1/x+1+x,1-x**2],x,2,4,10000): } \quad ,

getting a certain scheme with $801$ states.

For mod $32$, we type

{\tt gu:=CA([1/x+1+x,1-x**2],x,2,5,10000): } \quad ,

getting a certain scheme with $5093$ states.

For mod $64$, we type

{\tt gu:=CA([1/x+1+x,1-x**2],x,2,6,10000); } \quad ,

getting the output {\tt FAIL}, meaning that the number of needed states  exceeds our `cap', $10000$.

{\bf Fast Evaluation mod $p^a$}

Once an automatic $p$-scheme, $S$, is found for a combinatorial sequence modulo $p^a$, 
{\tt AutoSquared} can find very fast
the $N^{th}$ term of the sequence modulo $p^a$, for {\bf very large} $N$, using the 
procedure {\tt EvalLS(Z,N,i,p)}, with $i=1$.
For example, after first finding an  automatic $5$-scheme for the Motzkin numbers modulo $25$, by typing

{\tt gu:=CA([1/x+1+x,1-x**2],x,5,2,1000)[1]:} \quad ,

to get the remainder upon dividing $M_{10^{100}}$ by $25$, you should type:

{\tt EvalCA(gu,10**100,1,5);}

getting $12$. To get the first $N$ terms of the sequence (modulo $p^a$), once a scheme, {\tt S}, has been 
computed, type:

{\tt SeqCA(S,N,p)};

For example, with the above scheme (that we called {\tt gu}) (for the Motzkin numbers modulo $25$) 

{\tt SeqCA(gu,100000,5)};

takes {\tt 2.36} seconds to give you the first $100000$ terms,
and getting the first million terms, by typing ``{\tt SeqCA(gu,10**6,5)};'',
only takes {\tt 30} seconds.

\vfill\eject

{\bf  Congruence Linear Schemes}

The notion of {\it automatic $p$-scheme} defined above is conceptually attractive, 
since it can be modeled by a finite automaton with output. But, as can be
seen by the above example, the number of `states' (auxiliary functions) grows
very fast. But note that the space of polynomials modulo $p^a$ is a nice module
over the ring $\Z/(p^a \Z)$, and it is a shame to not take advantage of it.
So rather than waiting until no new pairs $[P(x),Q(x)]$ show up among the
``children'', it may be a good idea, whenever a new pair comes along, to see whether
it can be expressed as a {\bf linear combination} of previously encountered pairs with
the same $P(x)$ (which we already know stays the same after $a$ iterations, and only the
$Q(x)$'s change).

One can get away with many fewer auxiliary functions (`states') with the following notion.

{\bf Definition}: Let $\N$ be the set of non-negative integers, and let $p$ be a prime,
$a$ a positive integer, and let $M$ be a {\it module} over the ring of integers modulo $p^a$, $\Z/(p^a \Z)$.
A {\it linear $p$-scheme} for a function $f:\N \rightarrow M$ is
a set of {\it finitely} many (say $r$) auxiliary functions $A_1(n), \dots, A_r(n)$, where
$f(n)=A_1(n)$, and such that for each $i$ ($1 \leq i \leq r$), and each $\alpha$ ($0 \leq \alpha <p$), there exists a linear
combination
$$
A_i(pn+\alpha)= \sum_{j=1}^{r} C^{(\alpha)}_{i,j} A_{j}(n) \quad ,
$$
for some  $C^{(\alpha)}_{i,j} \in \{0,1, \dots, p^{a} -1 \}$, and there are initial conditions:
$$
A_i(0)=a_i \quad .
$$

Note that the previous notion of automatic $p$-scheme is the very special case, where for
each $\alpha$ and $i$,  there is exactly one $j$ (that equals $\sigma(\alpha,i)$) such that  $C^{(\alpha)}_{i,j}$ is non-zero,
and it has to be a $1$. 

{\bf Finding Linear $p$-Schemes in AutoSquared}

This is implemented, in {\tt AutoSquared}, by procedure {\tt LS} for single-variable $P$ and $Q$ and
by procedure {\tt LSmul} for multivariate $P$ and $Q$ (of course, {\tt LSmul} can handle also
a single variable, and we kept {\tt LS} both for old-time-sake and because it may be a bit faster for this special case).

The syntax for {\tt LS}  is

{\tt LS(Z,x,p,a,A,K): }

where {\tt Z} is a pair of single-variable functions $[P,Q]$, {\tt x} is the (single) variable 
name $x$ that serves as the argument of $P$ and $Q$, {\tt p} is a prime, {\tt a} is a positive integer, 
{\tt A} is a {\bf symbol} for expressing the linear expressions 
(where {\tt A[i]} means our humanese $A_i$),
and {\tt K} is (usually fairly large) positive integer,
stating the maximum number of `states' (auxiliary functions) that you are willing to put up with.
(It returns FAIL if the number of states exceeds {\tt K}.)

For example, to get a Linear  $2$-scheme for the Motzkin numbers, modulo $2$, 
(if you are willing to tolerate a scheme with up to $30$ members),

you type

{\tt gu:=LS([1/x+1+x,1-x**2],x,2,1,A,30): } \quad .

getting 

{\tt [[[A[2], A[2]], [A[3], A[4]], [A[3], A[3]], [0, A[2]]], [1, 1, 1, 0]]} \qquad ,

which is the {\bf same} as  the automatic $2$-scheme, spelled-out above, except it is phrased more verbosely.

If you type:

{\tt LS([1/x+1+x,1-x**2],x,2,2,A,30)[1];}

you would get the following linear $2$-scheme with $8$ states:

{\tt 
[[[A[2], A[8]], [A[3], A[7]], [A[4], A[5]], [A[4], A[6]], 

 [A[4], 2*A[3]+2*A[4]+3*A[5]], [3*A[4], 2*A[3]+2*A[4]+A[5]], [A[3]+A[4], A[2]+A[3]+A[4]], 

[A[3], 3*A[2]+A[3]+A[4]+3*A[5]]], 

[1, 1, 1, 1, 1, 3, 2, 1]]}  \qquad ,

that means that 
$$
A_1(2n)=A_2(n) \quad , \quad A_1(2n+1)=A_8(n) \quad ,  \dots \quad ,
$$
$$
A_8(2n)=A_3(n) \quad , \quad A_8(2n+1)=3A_2(n)+A_3(n)+A_4(n)+3A_5(n) \bmod 4 \quad .
$$
The corresponding automatic $2$-scheme has $24$ states.

For modulo $8$ we get $18$ states, compared to $128$ for the automatic $2$-scheme. For modulo $16$ we get $43$ states,
compared to $801$ states, and for modulo $32$ we get $96$ states, compared to $5093$ states.

Having  gotten a scheme, {\tt S}, phrased in terms of $A$, to get the first {\tt N} terms of the
sequence (modulo $p^a$), type

{\tt SeqLS(S,N,p,a,A) ;}

{\bf Other Highlights of AutoSquared}

Procedures {\tt BinCA}  and {\tt BinLS} find automatic $p$-schemes and linear $p$-schemes respectively
for {\it any} binomial coefficient sum. See the on-line help. 

As mentioned at the beginning, there are quite a few sample input and output files linked to from the front of this article

{\tt http://www.math.rutgers.edu/\~{}zeilberg/mamarim/mamarimhtml/meta.html} \quad \quad .

{\bf What about  congruences modulo integers that are NOT primes or prime powers?}

The Chinese Remainder Theorem comes to the rescue!

One first constructs as many automatic $p$-schemes, or linear $p$-schemes, for as many prime powers as
one could afford, or care about, and then one can very fast find the congruence class modulo
any integer involving these primes up to the given power.

{\bf The Maple packages {\tt CatalanLS}, {\tt MotzkinLS}, {\tt DelannoyLS}}

Using the main package {\tt AutoSquared}, our computer precomputed schemes for quite a few prime powers,
that enables us to find the remainder upon dividing by $m$, for many $m$, in particular, $m=1000$, getting
the last three digits of the Catalan, Motzkin, and Delannoy numbers given
at the prologue.

See the on-line  help in these packages.

{\bf Disclaimer}

Both the automatic $p$-schemes and the linear $p$-schemes that our Maple package output are
not guaranteed to be minimal. Of course the size does not change the fact that
they run in logarithmic time in the input, but the `implied constants' in the $O(\log n)$
algorithms are  most probably {\bf not} best-possible.

{\bf Conclusion}

The present project is {\bf yet another} {\it case study} in teaching computers to do
research all by themselves, once they were taught (programmed) the human tricks.
Once the computer mastered them, it can reproduce, in a few seconds, all the
previous results accomplished by humans, and go on to output much deeper results,
that no human, by himself, or herself, would  be able to do, hence getting,
much deeper results.
So the fact that the last three decimal digits of $M_{googol}$ are $187$, may not be as interesting
as Fermat's Last Theorem, but  is, {\it in some sense},  much deeper!

{\bf Acknowledgment}: The second-named author (DZ) was supported in part by a grant
from the National Science Foundation of the United States of America.
The authors thank Shalosh B. Ekhad for its many diligent and tedious computations and
proofs!

{\bf References}

[KKM] Manuel Kauers, Christian Krattenthaler,  and  Thomas W. M\"uller,
{\it A method for determining the mod-$2^k$ behaviour of recursive sequences, with applications to subgroup counting},
Electron.\ J. Combin.\ {\bf 18}(2) (2012), Article P37. \quad {\tt http://arxiv.org/abs/1107.2015}

[KM1] Christian Krattenthaler and  Thomas W. M\"uller,{\it A method for determining the mod-$3^k$ behaviour of recursive sequences},
preprint. \quad {\tt http://www.mat.univie.ac.at/\~{}kratt/artikel/3psl2z.html}

[KM2] Christian Krattenthaler and  Thomas W. M\"uller,
{\it A Riccati differential equation and free subgroup numbers for lifts of $PSL_2(\Z)$ modulo powers of primes},
J. Combin.\ Theory Ser.\ A {\bf 120} (2013), 2039--2063.
{\tt http://arxiv.org/abs/1211.2947}

[R] Gy\"orgy E. R\'ev\'esz, {\it ``Introduction to Formal Languages}, Dover, 1991.
[Originally published by McGraw-Hill, 1983].

[RY] Eric Rowland and Reem Yassawi, {\it Automatic congruences for diagonals of rational functions } \quad
{\tt http://arxiv.org/abs/1310.8635} \quad .

\bigskip

\hrule

\bigskip

Eric Rowland, Universit\'e du Qu\'ebec  \`a Montr\'eal,  Montr\'eal, Canada. \hfill\break
{\tt rowland at lacim dot ca}

Doron Zeilberger, Mathematics Department, Rutgers University (New Brunswick), Piscataway, NJ 08854, USA. \hfill\break
{\tt zeilberg at math dot rutgers dot edu}

{\bf Nov. 18, 2013}

\end